\documentclass[12pt,a4paper]{amsart}
\usepackage{amscd}
\usepackage{verbatim}
\usepackage[dvips]{graphicx}
\theoremstyle{plain}
  \newtheorem{thm}{Theorem}[section]
  \newtheorem{lem}[thm]{Lemma}

  \newtheorem{conj}[thm]{Conjecture}

\theoremstyle{definition}

\theoremstyle{remark}

\newcommand{\Z}{\mathbb{Z}}

\newcommand{\C}{\mathbb{C}}
\newcommand{\Q}{\mathbb{Q}}

\newcommand{\li}{\mathrm{Li}}
\newcommand{\I}{\sqrt{-1}}

\newcommand{\mA}{\mathcal A}
\newcommand{\mC}{\mathcal C}

\allowdisplaybreaks
\begin{document} 
\author
{Toshie Takata}
\address
[Toshie Takata]
{Department of Mathematics, Faculty of Science, Niigata University,
Niigata 950-2181, Japan}
\email
{takata@math.sc.niigata-u.ac.jp}
\title
{The colored Jones polynomial and the A-polynomial \\ 
 for twist knots}
\begin{abstract}
@We show that for a twist knot, the A-polynomial can be 
obtained from recurrences for the summand in Masbaum's formula 
of the colored Jones polynomial. 
Our result supports the {\it AJ conjecture} due to S.Garoufalidis. 
\end{abstract}
\maketitle

%%%%%%%%%%%%%%%%%%%%%%%%%%%%%%%%%%%%%%%%%%%%%%%%%%%%%%%%%%%%%%%%%%%%%
\section{introduction}
%%%%%%%%%%%%%%%%%%%%%%%%%%%%%%%%%%%%%%%%%%%%%%%%%%%%%%%%%%%%%%%%%%%%%

  The colored Jones polynomial is a one variable polynomial invariant 
 of a knot colored with irreducible finite dimensional 
 representation  of $sl(2,\C)$ \cite{Turaev:INVEM88}. 
  The original Jones polynomial \cite{Jones:link_polynomials}  
 corresponds to the colored Jones polynomial of a knot colored with  
 the irreducible 2-dimensional representation. 
 On the other hand, in \cite{Cooper/Culler/Gillet/Long/Shalen:INVEM94}, 
 D.Cooper, M.Culler, H.Gillet, D.D.Long, and P.R.Shalen introduced  
a  two variable polynomial invariant $A_K(l,m)$ of  a knot $K$, called 
the $A$-polynomial, by using the representations 
of the fundamental group of the complement of the knot into $SL_2(\mathbb{C})$.  
 In this paper, we prove that for a twist knot, the $A$-polynomial  
 can be  obtained from the colored Jones polynomial.
 
 Recently, 
 C.Frohman, R.Gelca, W.Lofaro \cite{Frohman/Gelca/Lofaro:TAMS01}
 introduced the peripheral ideal of a knot, and the noncommutative $A$-ideal 
 of a knot  as a generalization  of  the $A$-polynomial,
 via Kauffman bracket skein theory. A nontrivial element in the peripheral 
 ideal of a  knot  induces a recursive relation  of the Kauffman 
 bracket polynomial  of the knot, 
 which essentially equals to  the colored Jones polynomial of the  knot.
It is not yet proved that  for a knot 
except for the $(2,2p+1)$-torus knot and the figure eight knot, 
the peripheral ideal is nontrivial.

 Afterwards, in \cite{Garoufalidis/Le:math.GT}, 
S.Garoufalidis and T.T.Q.Le proved that the colored Jones polynomial  of a knot
satisfies a nontrivial recursive relation. 
 Moreover, in \cite{Garoufalidis:math.GT}, 
 Garoufalidis defined the recursion ideal,  which is 
identified with the set of  recursive relations of the colored 
Jones polynomial, and defined the noncommutative 
$A$-polynomial $A_q(K)$ of a knot $K$ called $A_q$-{\it polynomial},  as 
a  generator of the ideal. 
Let us denote by $J_K(n)$ the $n$-colored Jones polynomial associated 
  with  irreducible $n$-dimensional representation of  $sl(2,\C)$. 
We consider two operators $E$ and $Q$ acting on $J_K(n)$ 
by $EJ_K(n)=J_K(n+1)$ and $QJ_K(n)=q^nJ_K(n)$. 
The element $A_q(K)$ is of the form  
$A_q(K)(E,Q)=\sum_k a_k E^k$ with $a_k$ in $\Z [q,Q]$.
Then, Garoufalidis conjectured 
\begin{conj} [{\it The AJ conjecture}] 
For every knot $K$ in $S^3$, $A_K(l,m)=\varepsilon A_q(K)(l,m^2)$ , 
where $\varepsilon $ is the evaluation map at $q=1$. 
\end{conj}
Furthermore, he showed that the conjecture holds for  
 the trefoil knot and the figure eight knot,  
 by using the mathematica package {\tt qZeil.m} develped by  Paul and Riese 
 \cite{Paule/Riese:FIC97} \cite{Paule/Riese:qZeil}, to find 
  a recursive relation of the colored Jones polynomial. 
 
  We focus on twist knots, for which we have 
 the formula of the colored Jones polynomial obtained by G.Masbaum 
  in \cite{Masbaum:AGT03} and the formula of the $A$-polynomial 
  obtained by J.Hoste and P.D.Shanahan in \cite{Hoste/Shanahan:JK}.  
  The purpose of this paper is to prove that for a twist knot, 
   the $A$-polynomial can be obtained from recursive relations of 
  the summand in Masbaum's formula of the colored Jones polynomial. 
  Our result supports   the AJ conjecture for a twist knot.
  Moreover, using the mathematica package 
  {\tt qMultisum.m} developed by Reise \cite{Riese:JSC03}, 
   we observe that the $A$-polynomial for the knots $5_2$ and $6_1$ 
  can be obtained from a recursive relation of the colored Jones 
  polynomial. 
  We will also discuss a relation between our result and 
  `Volume conjecture' due to R.M.Kashaev, H. Murakami and J. Murakami
   (\cite{Kashaev:LETMP97},\cite{Murakami/Murakami:volume}).

  This paper is organized as follows. 
  In Section 2, we recall  Masbaum's formula of the colored Jones polynomial 
  for a twist knot and state Main theorem (Theorem \ref{thm:main}). In Section 3, we prove Theorem \ref{thm:main}.
   In Section 4, we relate  our result 
  to  the AJ conjecture and give a result about recursive relations 
  of the Jones polynomials for the knots $5_2$ and $6_1$.

 The author would like to thank Y. Yokota for 
  helpful conversation. 
  
%%%%%%%%%%%%%%%%%%%%%%%%%%%%%%%%%%%%%%%%%%%%%%%%%%%%%%%%%%%%%%%%%
\section{Main Theorem}
%%%%%%%%%%%%%%%%%%%%%%%%%%%%%%%%%%%%%%%%%%%%%%%%%%%%%%%%%%%%%%%%%

 We start with the review of the definition of the $A$-polynomial of a knot 
$K$ in $S^3$ (see \cite{Cooper/Culler/Gillet/Long/Shalen:INVEM94} ). 
Let  $M_K$ be the complement of $K$ and $R(M_K)$ the set 
of all representations of $\pi_1(M_K)$ into $SL(2,\C)$. 
The set $R(M_K)$ is an affine algebraic variety. 
Noting that $SL(2,\C)$ acts on representations by conjugation, 
we restrict our attention to the subset $R_U$ of $R(M_K)$ consisting of 
a represetation $\rho$ satisfying that $\rho(\mu)$ and $\rho(\lambda)$ 
are upper triangular matirices 
for the meridian  $\mu$ and the preferred longitude $\lambda$ of $K$.
We can define a projection $\xi$ from $R_U$ to $\C^2$ by 
$\xi(\rho)=(l,m)$ 
for $\rho \in R(M_K)$ with 
$$ \rho(\mu)= \begin{pmatrix}
              m & * \\
              0  & m^{-1}
            \end{pmatrix}, \; 
 \rho(\lambda)= \begin{pmatrix}
              l & * \\
              0 & l^{-1} \\
            \end{pmatrix} .
$$  
The Zariski closure of $\xi(R_U)$ has a structure of an algebraic variety
 in $\C^2$ and each of its irreducible complex-dimension-one components 
 is a curve, which is defined by a polynomial with integer coefficients in $l$ and $m$. 
The $A$-polynomial is the product of those defining polynomials. 
We note that the $A$-polynomial of $K$ has a factor $l-1$, which correponds to 
abelian representations. 
So, we denote by $A_K(l,m)$ the $A$-polynomial divided by $l-1$.

 Nextly, we recall Masbaum's formula of the colored Jones polynomial for 
 a twist knot. 
Some notations are fixed.
\begin{eqnarray*}
& & \{n\}=s^{n}-s^{-n},  \; s^2=q, \\ %[n]=\frac {\{n\}}{\{1\}}, 
& & \{n\}!=\{n\}\{n-1\}\cdots \{1\}, \; (q)_n=(1-q)(1-q^2)\cdots (1-q^n).
\end{eqnarray*}
Let  $J_K(n)$ be the colored Jones polynomial of a 0-framing knot 
$K$ colored with the $n$-dimensional irreducible representation 
of $sl_2(\C)$, where $J_K(n)$ is  normalized by $J_{\bigcirc}(n)=1$ 
for the trivial 
knot ${\bigcirc}$. 
($J'_K(n)$ in \cite{Masbaum:AGT03} is equal to $J_K(n)$ in this paper.)

 Let $K_p$ be the twist knot drawn in figure 1.
In \cite{Masbaum:AGT03}, Masbaum obtained the following formula for the colored Jones polynomial of the twist knot $K_p$. 
\begin{thm}\cite{Masbaum:AGT03}
The colored Jones polynomial of the twist knot $K_p$ is given by 
\begin{equation}
J_{K_p}(n)=\sum_{k=0}^{\infty} \mC_{K_p}(k) 
       \frac {\{n-k \} \{n-k+1\} \cdots \{n+k\}} {\{n\}},
\label{eqn:mas}
\end{equation}
where
\begin{equation}
\mC_{K_p}(k)
=(-1)^{k+1}s^{k(k+3)/2}\sum_{l=0}^k (-1)^l q^{l(l+1)p} \{2l+1\}
 \frac {\{k \}! }{\{k+l+1 \}! \{k-l\}!}. 
\end{equation}
\end{thm}

 We make use of the rearranged Masbaum's formula 
 (\cite{Garoufalidis/Le:math.GT}) 
\begin{eqnarray*}
J_{K_p}(n)&=&\sum_{k=0}^{\infty}\sum_{l=0}^k (-1)^{k+1} q^{k(k+3)/2}
        q^{nk}\frac {(q^{-1})_{n+k}(q^{-1})_{n-1}}
                    {(q^{-1})_{n}(q^{-1})_{n-k-1}}\\
 {} & & \qquad \times 
       (-1)^l q^{l(l+1)p+l(l-1)/2}(q^{2l+1}-1) 
       \frac {(q)_k}{(q)_{k+l+1} (q)_{k-l}}.\qquad \qquad
\end{eqnarray*}
Let us put 
\begin{eqnarray*}
& & F(n,k,l)\\
& &=(-1)^{k+1} q^{k(k+3)/2}
        q^{nk}\frac {(q^{-1})_{n+k}(q^{-1})_{n-1}}
                    {(q^{-1})_{n}(q^{-1})_{n-k-1}} 
       (-1)^l q^{l(l+1)p+l(l-1)/2}(q^{2l+1}-1) 
       \frac {(q)_k}{(q)_{k+l+1} (q)_{k-l}}. 
\end{eqnarray*}
Then, it easily follows that 
\begin{eqnarray}
f_0(q,q^n,q^k,q^l) & := &\frac {F(n+1,k,l)} {F(n,k,l)}
 =q^k\frac{(1-q^{-n-k-1})(1-q^{-n})}{(1-q^{-n-1})(1-q^{-n+k})},
 \label{eqn:f0}\\
f_1(q,q^n,q^k,q^l) & := &\frac {F(n,k+1,l)} {F(n,k,l)}
 =-q^{k+n+2}
  \frac {(1-q^{-n-k-1})(1-q^{-n+k+1})(1-q^k)}{(1-q^{k+l+2})(1-q^{k-l+1})},
 \label{eqn:f1} \\ 
f_2(q,q^n,q^k,q^l) & := &\frac {F(n,k,l+1)} {F(n,k,l)}
  =-q^{(2p+1)l+2p}
  \frac {(q^{2(l+1)+1}-1)(1-q^{k-l})}{(q^{2l+1}-1)(1-q^{k+l+2})}.
 \label{eqn:f2}
\end{eqnarray}  
 We  obtain the following result.
\begin{thm}\label{thm:main}
From the three equations 
\[
f_0(1,m^2,x,y)=l, \; f_1(1,m^2,x,y)=1, \; f_2(1,m^2,x,y)=1,
\] 
we can obtain 
$$x=\frac {lm^2+1}{m^2+l}, \quad 
(y+1)(1-x)y^p h_p(m,x)=0,  
$$
where $h_p(m,x)$ is an element in $\Z [m^{\pm 1},x]$ and satisfies that 
if we put 
$$B_{K_p}(l,m):=\left\{
    \begin{array}{ll}
    (m^2+l)^{2p-1}(m^2)^p h_p(m,\frac {lm^2+1}{m^2+l} ), \quad & p>0,\\
    (m^2+l)^{2|p|}(m^2)^{|p|} h_p(m,\frac {lm^2+1}{m^2+l} ),  & p \le 0,\\
    \end{array}
    \right.
 $$
then, $B_{K_p}(l,m)=A_{K_p}(l,m)$.    
\end{thm}

%%%%%%%%%%%%%%%%%%%%%%%%%%%%%%%%%%%%%%%%%%%%%%%%%%%%%%%%%%%%%%%%%%%
\section{Proof of Theorem \ref{thm:main}}
%%%%%%%%%%%%%%%%%%%%%%%%%%%%%%%%%%%%%%%%%%%%%%%%%%%%%%%%%%%%%%%%%

 Let $K_p$ be the twist knot pictured in Figure 1.
 To prove Theorem \ref{thm:main}, we use the following  formula for       
the $A$-polynomial of the twist knot $K_p$ due to 
Hoste and Shanahan in \cite{Hoste/Shanahan:JK}.

\begin{thm}(Hoste-Shanahan \cite{Hoste/Shanahan:JK})
For $p\ne -1,0,1,2$, $A_{K_p}(l,m)$ is given recusively by 
\begin{equation}
A_{K_p}(l,m)=\left\{ 
            \begin{array}{ll}
            c A_{K_{p-1}}(l,m)-d A_{K_{p-2}}(l,m), & p>0,\\
            c A_{K_{p+1}}(l,m)-d A_{K_{p+2}}(l,m), & p<0, 
            \end{array}
            \right.\label{eqn:HS}
\end{equation}
where 
\begin{eqnarray*}
c&=&-l+l^2+2lm^2+m^4+2lm^4+l^2m^4+2lm^6+m^8-lm^8, \\
d&=&m^4(l+m^2)^4,
\end{eqnarray*}
and with initial conditions 
\begin{eqnarray*}
A_{K_2}(l,m)&=&-l^2+l^3+2l^2m^2+lm^4+2l^2m^4-lm^6-l^2m^8 \\
 {} & &        +2lm^{10}+l^2m^{10}+2lm^{12}+m^{14}-lm^{14},\\
A_{K_1}(l,m)&=  & l+m^6,\\
A_{K_0}(l,m)&=  & 1, \\
A_{K_{-1}}(l,m)&=  &-l+lm^2+m^4+2lm^4+l^2m^4+lm^6-lm^8.\\
\end{eqnarray*}
\end{thm}

\begin{figure}
\centering
\includegraphics[scale=0.5]{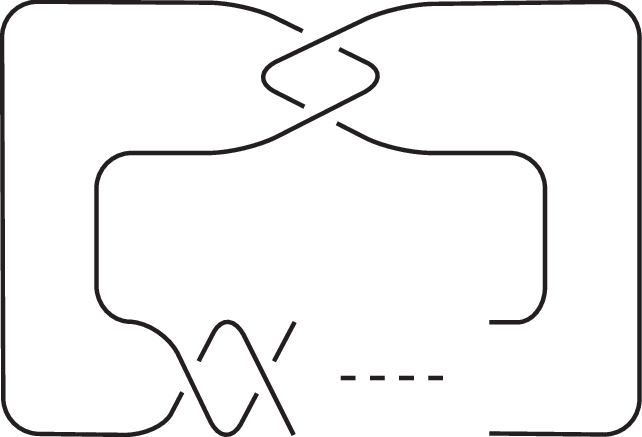}

p full twists
\begin{center}
 Figure 1
\end{center}
\end{figure} 

Now, we give a proof of Theorem \ref{thm:main}. 

\quad

\noindent 
{\it Proof of Theorem \ref{thm:main}}
From the first equation 
\[
f_0(1,m^2,x,y)=x\frac {(1-\frac 1 {xm^2})}{(1-\frac x {m^2})}=l, 
\]
we get 
\[
x=\frac {lm^2+1}{m^2+l}.
\]
Moreover,  the second equation 
\[
f_1(1,m^2,x,y)=-xm^2
               \frac {(1-\frac 1 {xm^2})(1-\frac x {m^2})(1-x)}
                     {(1-xy)(1-\frac x y)}
                     =1
\]
gives  
\begin{equation}
y^2+1=\frac y {m^2}(m^4-xm^4+x^2m^2+m^2+1-x). \label{eqn:y}
\end{equation}
It is clear that the claim for $p=0$ holds. 
We consider the case $p>0$. 
The third equation
\[f_2(1,m^2,x,y)=-y^{2p+1}\frac {(1-\frac x y)}{(1-xy)}=1 
 \]
implies that 
\begin{equation}
y^{2p+1}+1-y^{2p}x-xy=0, 
\end{equation}
which can be changed to 
\begin{eqnarray*}
0&=&(y+1)(y^{2p}-y^{2p-1}+y^{2p-2}- \cdots +1-xy(y^{2p-2}-y^{2p-3}+\cdots +1))\\
&=&(y+1)\{(y^2+1-(1+x)y)(y^{2(p-1)}+y^{2(p-2)}+\cdots+y^2+1) \\
& &\hspace{2truecm} -(1-x)(y^{2(p-1)}+y^{2(p-2)}+\cdots+y^2)\}.
\end{eqnarray*}
From  the equation (\ref{eqn:y}), we can get 
\begin{equation}
(y+1)(1-x)H_p(m,x,y)=0, 
\end{equation}
where $H_0(m,x,y)=1$ and for $p \ge 1$, 
\begin{eqnarray*}
& &H_p(m,x,y)\\
& &=\frac y {m^2}(m^4-xm^2+1)(y^{2(p-1)}+y^{2(p-2)}+\cdots+y^2+1)
              -(y^{2(p-1)}+y^{2(p-2)}+\cdots+y^2).
\end{eqnarray*}

To prove the theorem for $p>  0$, we need 

\begin{lem}\label{lem:key}
We have $H_p(m,x,y)=y^p h_p(m,x)$, where 
for $p\ge 2$, $h_p(m,x)$ is given recursively by
\begin{equation}
h_p(m,x)=\frac a {m^2} h_{p-1}(m,x)-h_{p-2}(m,x) 
\end{equation}
with $a=m^4-xm^4+x^2m^2+m^2+1-x $, 
$h_0(m,x)=1$, $h_1(m,x)=\frac b {m^2}$, and $b=m^4-xm^2+1 $.
\end{lem}

\noindent 
{\it Proof of lemma \ref{lem:key}}
We show by induction on $p$. 
When $p=1$, it follows that 
$H_1(m,x,y)=\frac y {m^2}(m^4-xm^2+1)=y h_1(m,x)$. 
In the case $p=2$, we obtain that  
\begin{eqnarray*}
H_{2}(m,x,y)& =& H_1(m,x,y)+\frac {yb} {m^2} y^2-y^2 \\
              & =& yh_1(m,x)+y^3 h_1(m,x)-y^2\\
              & =& y (h_1(m,x)(y^2+1)-y)\\
              & =& y (\frac {ay} {m^2}h_1(m,x)-y)\\
              & =& y^2 (\frac {a} {m^2}h_1(m,x)-1)\\
              & =& y^2 (\frac {a} {m^2}h_1(m,x)-h_0(m,x)).\\
\end{eqnarray*}
Here we used the equation (\ref{eqn:y}) in the fourth equality. 
Suppose that it holds for $p \le n$. 
By hypothesis of induction and the definition of $H_p(m,x,y)$, we can 
obtain  
\begin{eqnarray*}
H_{n+1}(m,x,y)& =& H_n(m,x,y)+\frac {yb} {m^2} y^{2n}-y^{2n}\\
              & =& y^n h_n(m,x,y)+y^2(\frac {yb} {m^2} y^{2(n-1)}-y^{2(n-1)})\\
              & =& y^n h_n(m,x,y)+y^2(H_n(m,x,y)-H_{n-1}(m,x,y))\\
              & =& y^n h_n(m,x,y)+y^2(y^n (h_n(m,x,y)-y^{n-1} h_{n-1}(m,x,y))\\
              &=& y^n((y^2+1) h_n(m,x,y)-y h_{n-1}(m,x,y))\\
              &=& y^n(\frac {a} {m^2}y h_n(m,x,y)-y h_{n-1}(m,x,y))\\
              &=& y^{n+1}(\frac {a} {m^2} h_n(m,x,y)- h_{n-1}(m,x,y)).
\end{eqnarray*}
This completes the proof. \qed

Let us  go back to the proof of the theorem. 
From  Lemma \ref{lem:key}, it is clear that $h_p(m,x)$ 
is in $\Z [m^{\pm 1},x]$. 
We put 
\begin{equation}
B_{K_p}(l,m):=(m^2+l)^{2p-1}(m^2)^p h_p(m,\frac {lm^2+1}{m^2+l} ),
\end{equation}
and then, it is easy to see that 
$B_{K_p}(l,m)$ satisfies the recursive relation 
\begin{equation}
B_{K_p}(l,m)=a(m^2+l)^2 B_{K_{p-1}}(l,m)-m^4(m^2+l)^4B_{K_{p-2}}(l,m).
\end{equation}
By comparing it 
with the recursive relation (\ref{eqn:HS}) of $A_{K_p}(l,m)$ 
 by Hoste and Shanahan, 
and by noting that 
$a(m^2+l)^2=-l+l^2+2lm^2+m^4+2lm^4+l^2m^4+2lm^6+m^8-lm^8=c$ and 
$d=m^4(m^2+l)^4$,
it is shown  that $B_{K_p}(l,m)=A_{K_p}(l,m)$ for $p\ge 0$. 
In a similar way to the case $p \ge 0$, 
we can prove the theorem in the case $p<0$. \qed
%\end{proof}

%%%%%%%%%%%%%%%%%%%%%%%%%%%%%%%%%%%%%%%%%%%%%%%%%%%%%%%%%%%%%%%%%
\section{Observations}
%%%%%%%%%%%%%%%%%%%%%%%%%%%%%%%%%%%%%%%%%%%%%%%%%%%%%%%%%%%%%%%%%
 
 We start with reviewing  Wilf-Zeilberger's algorithm. 
  For more details, 
we refer to \cite{Garoufalidis/Le:math.GT} and \cite{Garoufalidis:math.GT}.
 A disctrete function $F(n,{\bf k})$ is called proper $q$-hypergeometric 
 if it is one of the form 
\[
 F(n,{\bf k})=\frac {\prod_{s}(A_s;q)_{a_s n+{\bf b}_s\cdot {\bf k} +c_s}}
                    {\prod_{t}(B_t;q)_{u_t n+{\bf v}_t\cdot {\bf k} +w_t}}
                    q^{A(n,{\bf k})}\xi^{\bf k},
\]
where ${\bf k}=(k_1,\ldots,k_r)$, 
$A_s, B_t \in \Q (q)$, $a_s$, $u_t$ are integers, 
${\bf b}_s$, ${\bf k}_s$ are vectors of $r$ integers, 
$A(n,{\bf k})$ is a quadratic form, $c_s$, $w_t$ are variables 
and $\xi$ is an $r$ vector of elements in $\Q (q)$, and for $k\ge 0$, 
\[
(A,q)_k:=\left\{
        \begin{array}{ll}
        (1-A)(1-Aq)\cdots (1-Aq^{k-1}) & \mathrm{if} \; k>0,\\
        1                              & \mathrm{if} \; k=0.\\
        \end{array}
        \right.
\]

%  we will comment on a relation between our result and the AJ conjecture 
%  due to Garoufalidis \cite{Garoufalidis:math.GT}. 
 In \cite{Garoufalidis/Le:math.GT}, Garoufalidis and Le proved that 
the colored Jones polynomial $J_K(n)$ can be written as a multisum 
\[
J_K(n)=\sum_{\bf k=0}^{\infty} F(n,k_1,\ldots,k_r)
\]
of  a proper $q$-hypergeometric function $F(n,k_1,\ldots,k_r)$, 
where only finitely many terms 
in the right hand side are nonzero.
Consider the operators $E$, $E_i$, $Q$ and $Q_j$, $(1 \le i,j \le r)$  
 acting on a discrete function $f: \mathbb{N}^r \to \Z[q^{\pm }]$ by 
\begin{eqnarray*}
 & &(E f)(n,k_1,\ldots,k_r)=f(n+1,k_1,\ldots,k_r), \\
 & &(E_i f)(n,k_1,\ldots,k_r)=f(n,k_1,\ldots,k_{i-1},k_i+1,k_{i+1},\ldots,k_r), \\
 & &(Q f)(n,k_1,\ldots,k_r)=q^n f(n,k_1,\ldots,k_r),\\ 
 & &(Q_j f)(n,k_1,\ldots,k_r)=q^{k_j} f(n,k_1,\ldots,k_r). 
\end{eqnarray*}
Then, they satisfy the  relations
\begin{eqnarray*}
& & Q Q_i=Q_i Q,  Q_i Q_j=Q_j Q_i, E E_i=E_i E, E_i E_j=E_j E_i,\\
& & Q_i E=E Q_i,  Q_i E_j= E_j Q_i \; \mathrm{for} \; i\ne j, 
 E Q =q Q E,  E_i Q_i =q Q_i E_i, 
\end{eqnarray*}
%\end{comment}
which are denoted by $(Rel_q)$. 
The $q$-Weyl algebra $\mA$ is defined to be a noncomutative algebra 
with presentation 
\[
\mA=\frac {\Z [q^{\pm}]\langle Q, Q_1,\ldots,Q_r,E,E_1,\ldots,E_r \rangle}
           {(Rel_q)}.
\]

Wilf and Zeilberger    proved that 
\begin{thm}\cite{Wilf/Zeilberger:INVEM92}
Every proper $q$-hypergeometric function $F(n,{\bf k})$ satisfies 
a ${\bf k}$-free recurrence 
\begin{equation}
\sum_{i,{\bf j} \in S} \sigma_{i,{\bf j}}(q^n) F(n+i,{\bf k}+{\bf j})=0,
\label{eqn:k-free}
\end{equation}
where $S$ is a finite set, ${\bf k}=(k_1,\ldots,k_r)$, 
${\bf j}=(j_1,\ldots,j_r)$, and $\sigma_{i,{\bf j}}(q^n)$ are 
polynomials in $q^n$ 
with coefficients in $\Q (q)$. 
\end{thm}
Putting 
\[
P=P(E,Q,E_1,\ldots, E_r)
=\sum_{i,{\bf j} \in S} \sigma_{i,{\bf j}}(Q)E^i E^{\bf j},
\]
where $E^{\bf j}={E_1}^{j_1}\cdots {E_r}^{j_r}$, 
the recurrence (\ref{eqn:k-free}) can be written in operator notation as $PF=0$. 

 Moreover, they showed that any ${\bf k}$-free recurrence 
 can be transformed to 
\[
P=P(E,Q,1,\ldots,1)+\sum_{i=1}^r (E_i-1) R_i(E,Q,E_1, \ldots,E_r).
\]
Applying  it to $F(n,k_1,\ldots,k_r)$ and summing over ${\bf k} \ge {\bf 0}$ 
give us 
\[
P(E,Q,1,\ldots,1)J_K(n)=error(n), 
\]
where $error(n)$ is a sum of a multisum of proper $q$-hypergeometric function 
with $r-1$ variables. So, repeating the process, 
we arrive at  a homogeneous recursive relation of $J_K(n)$
\[
P_1(E,Q)P(E,Q,1,\ldots,1)J_K(n)=0,
\]
where $P_1(E,Q)$ is a polynomial in $E$ with coefficients in $\Q (q) [Q]$.

On the other hand, 
from the fact that $F(n,k_1,\ldots,k_r)$ is proper $q$-hypergeometric, 
we can write $E F/F=A/B$ and $E_i F/F=A_i/B_i$ for polynomials 
$A, B, A_i, B_i \in \Q (q) [q^n,q^{k_1},\ldots,q^{k_r}]$. 
Then, putting $Q=q^n$ and $Q_i=q^{k_i}$, it follows 
that $(BE-A)F=0$, $(B_iE_i-A_i)F=0$ and 
that $BE-A$, $B_iE_i-A_i$ generate the annihilation ideal of $F$ in $\mA$. 
From these generators, the Wilf-Zeilberger's algorithm finds a 
${\bf k}$-free recurrence of $F(n,{\bf k})$ and then 
 a  recursive relation of the sum of $\sum_{\bf k} F(n,{\bf k})$. 

Now, we apply the above algorithm to $J_{K_p}(n)$. 
From the equations (\ref{eqn:f0}), (\ref{eqn:f1}) and (\ref{eqn:f2}), 
setting 
\begin{eqnarray*}
 & B=(1-q^{-n-1})(1-q^{-n+k}), & \;A=q^k(1-q^{-n-k-1})(1-q^{-n}),\\
 & B_1=(1-q^{k+l+2})(1-q^{k-l+1}), 
 & A_1=-q^{k+n+2}(1-q^{-n-k-1})(1-q^{-n+k+1})(1-q^k),\\
 & B_2=(q^{2l+1}-1)(1-q^{k+l+2}),
 & A_2=q^{(2p+1)l+2p}(q^{2(l+1)+1}-1)(1-q^{k-l}),
\end{eqnarray*}
we have $(BE-A)F=0$, $(B_i E_i-A_i)F=0$ for $i=1,2$, with $Q=q^n$, 
$Q_1=q^k$ and $Q_2=q^l$.
Identifying $Q$ with $m^2$, $E$ with $l$, $Q_1$ with $x$, $Q_2$ with $y$, 
from the definition of $f_i$ for $i=0,1,2$ in Section 2,
we observe that $\varepsilon (BE-A)$ corresponds to 
 $f_0(1,m^2,x,y)-l$ and that $\varepsilon (B_iE_i-A_i)$ with $E_i=1$ 
 corresponds to $ f_i(1,m^2,x,y)-1$ for $i=1,2$,  
 where $\varepsilon$ is the evaluation map at $q=1$ 
 defined in \cite{Garoufalidis:math.GT}.
 So, our result may support the  AJ conjecture.

Nextly, we will present the usage of the mathematica package {\tt qMultiSum.m}
 with the above theorical background,  which is  
 developed by A. Riese \cite{Riese:JSC03}, to obtain inhomogeneous recurrences    of the colored Jones polynomial  of the knots $5_2$ and $6_1$. 
 We will also relate them  
 to the $A$-polynomial.
 
%%%%%%%%%%%%%
% 5_2
%%%%%%%%%%%%

 The formula in \cite{Masbaum:AGT03} allows us to write 
\begin{equation}
J_{5_2}(n)=\sum_{k=0}^{\infty}\sum_{l=0}^k (-1)^{k+1} q^{(3k^2+5k)/2}
       q^{nk}\frac {(q^{-1})_{n+k}(q^{-1})_{n-1}}
                    {(q^{-1})_{n}(q^{-1})_{n-k-1}}
       q^{-l(k+1)}\frac {(q^{-1})_k}{(q^{-1})_{l} (q^{-1})_{k-l}}. 
\end{equation}
We denote  the summand by $F(n,k,l)$. 
Then, the mathematica program soft  
{\tt qMultiSum.m} computes the ${\bf k}$-free recurrence of $F(n,k,l)$
\begin{eqnarray*}
& & -q^{4 + 7n}\left( q - q^n \right)\left( q^2 - q^n \right) 
       \left( q^5 - q^n \right) \left( q + q^n \right) 
       \left( q^2 + q^n \right)        \left( q - q^{2n} \right) 
       \left( q^3 - q^{2n} \right) F(-5 + n,k,l) \\
& &   + q^{1 + 3n}\left( q - q^n \right) \left( q^2 - q^n \right) 
      {\left( q^4 - q^n \right) }^2
     \left( q + q^n \right) \left( q^2 + q^n \right) \left( q^4 + q^n \right)
      \left( q - q^{2n} \right) \left( q^3 - q^{2n} \right) \\
& &  \times \left( q^9 - q^{2n} \right) 
     \left( q^4 - q^n - q^{1 + n} \right) F(-4 + n,-1 + k,-1 + l) \\
& & - q^{9 + 2n}\left( q - q^n \right) \left( q^2 - q^n \right) 
    {\left( q^4 - q^n \right) }^2
     \left( q + q^n \right) \left( q^2 + q^n \right) \left( q^4 + q^n \right) 
     \left( q - q^{2n} \right) \\
& &  \times \left( q^3 - q^{2n} \right) 
     \left( q^9 - q^{2n} \right)  F(-4 + n,-1 + k,l)\\
& & - q^{5 + 6n}\left( 1 + q + q^2 + q^3 + q^4 \right) 
     \left( q - q^n \right) {\left( q^4 - q^n \right) }^2\left( q + q^n \right)      \left( q^4 + q^n \right) \left( q - q^{2n} \right) \\
& &  \times \left( q^3 - q^{2n} \right)  F(-4 + n,k,l)\\
& & + \left( q - q^n \right) \left( q^2 - q^n \right) 
     {\left( q^3 - q^n \right) }^3\left( q^4 - q^n \right) 
     \left( q + q^n \right) 
     \left( q^2 + q^n \right) \left( q^3 + q^n \right) \left( q^4 + q^n \right)      \left( q - q^{2n} \right) \\
& & \times \left( q^3 - q^{2n} \right) 
     \left( q^7 - q^{2n} \right) 
     \left( q^9 - q^{2n} \right) F(-3 + n,-2 + k,-2 + l)\\
& & + q^{2 + 2n}\left( q - q^n \right) {\left( q^3 - q^n \right) }^2
     \left( q^4 - q^n \right) 
     \left( q + q^n \right) \left( q^3 + q^n \right) \left( q^4 + q^n \right) 
     \left( q - q^{2n} \right) \left( q^3 - q^{2n} \right)\\
& &  \times \left( q^9 - q^{2n} \right) 
     \left( q^5 + q^6 + q^{2n} - q^{1 + n} - 2q^{2 + n} - 2q^{3 + n} - q^{4 + n}    \right) F(-3 + n,-1 + k,-1 + l)\\ 
& & - q^{5 + n}\left( q - q^n \right) {\left( q^3 - q^n \right) }^2
     \left( q^4 - q^n \right) \left( q + q^n \right) \left( q^3 + q^n \right) 
     \left( q^4 + q^n \right) \left( q - q^{2n} \right) 
     \left( q^3 - q^{2n} \right) \\
& &  \times 
     \left( q^9 - q^{2n} \right) \left( q^6 + q^{2n} + q^{1 + 2n} \right) 
     F(-3 + n,-1 + k,l)\\
& & - q^{6 + 5n}\left( 1 + q^2 \right) \left( 1 + q + q^2 + q^3 + q^4 \right) 
     \left( q - q^n \right) {\left( q^3 - q^n \right) }^2\left( q + q^n \right)      \left( q^3 + q^n \right) \left( q - q^{2n} \right)\\
& &  \times \left( q^9 - q^{2n} \right) F(-3 + n,k,l)\\
& & - q\left( q - q^n \right) {\left( q^2 - q^n \right) }^3
     \left( q^3 - q^n \right) \left( q^4 - q^n \right) \left( q + q^n \right) 
     \left( q^2 + q^n \right) \left( q^3 + q^n \right) \left( q^4 + q^n \right)      \\
& & \times  \left( q - q^{2n} \right)\left( q^3 - q^{2n} \right) 
     \left( q^7 - q^{2n} \right) 
     \left( q^9 - q^{2n} \right) F(-2 + n,-2 + k,-2 + l)\\
& & + q^{4 + n}\left( q - q^n \right) {\left( q^2 - q^n \right) }^2
    \left( q^4 - q^n \right) 
     \left( q + q^n \right) \left( q^2 + q^n \right) \left( q^4 + q^n \right) 
     \left( q - q^{2n} \right) \left( q^7 - q^{2n} \right)
      \\
& &  \times  \left( q^9 - q^{2n} \right)
     \left( q^5 + q^{2n} - q^{1 + n} - 2q^{2 + n} - 2q^{3 + n} - q^{4 + n}
      + q^{1 + 2n} \right) F(-2 + n,-1 + k,-1 + l)\\
& & - q^{3 + 2n}\left( q - q^n \right) {\left( q^2 - q^n \right) }^2
    \left( q^4 - q^n \right) \left( q + q^n \right) \left( q^2 + q^n \right)
     \left( q^4 + q^n \right) 
     \left( q - q^{2n} \right) \left( q^7 - q^{2n} \right)\\
& &  \times \left( q^9 - q^{2n} \right) 
     \left( q^4 + q^5 + q^{2n} \right) F(-2 + n,-1 + k,l)\\
& & - q^{7 + 4n}\left( 1 + q^2 \right) \left( 1 + q + q^2 + q^3 + q^4 \right) 
     {\left( q^2 - q^n \right) }^2\left( q^4 - q^n \right) 
     \left( q^2 + q^n \right) 
     \left( q^4 + q^n \right) \left( q - q^{2n} \right) \\
& &  \times \left( q^9 - q^{2n} \right)  F(-2 + n,k,l)\\
& & - q^{6 + n}{\left( q - q^n \right) }^2\left( q + q^2 - q^n \right) 
     \left( q^3 - q^n \right) \left( q^4 - q^n \right) \left( q + q^n \right) 
     \left( q^3 + q^n \right) \left( q^4 + q^n \right) 
     \left( q - q^{2n} \right) \\
& &  \times \left( q^7 - q^{2n} \right) \left( q^9 - q^{2n} \right) 
     F(-1 + n,-1 + k,-1 + l) \\
& &- q^{5 + 3n}{\left( q - q^n \right) }^2\left( q^3 - q^n \right) 
     \left( q^4 - q^n \right) 
     \left( q + q^n \right) \left( q^3 + q^n \right) \left( q^4 + q^n \right) 
     \left( q - q^{2n} \right) \left( q^7 - q^{2n} \right)\\
& &  \times \left( q^9 - q^{2n} \right)  F(-1 + n,-1 + k,l)\\
& & - q^{8 + 3n}\left( 1 + q + q^2 + q^3 + q^4 \right) 
     {\left( q - q^n \right) }^2\left( q^4 - q^n \right) \left( q + q^n \right)      \left( q^4 + q^n \right) \left( q^7 - q^{2n} \right) \\
& &  \times \left( q^9 - q^{2n} \right) F(-1 + n,k,l)\\
& & + q^{9 + 2n}\left( q^3 - q^n \right) \left( q^4 - q^n \right) 
     \left( -1 + q^n \right) \left( q^3 + q^n \right) \left( q^4 + q^n \right) 
     \left( q^7 - q^{2n} \right) \\
& & \times \left( q^9 - q^{2n} \right) F(n,k,l) = 0.
\end{eqnarray*}
Moreover, it can be converted into the inhomogeneous recursive relation  
of the colored Jones polynomial of the knot $5_2$ 
\begin{eqnarray*}
& &q^{9 + 7n} J_{5_2}(n)\\
& &+ \frac {{\left( -1 + q^{1 + n} \right) }^2\left( 1 + q^{1 + n} \right)}
           {\left( -1 + q^n \right) \left( -1 + q^{3 + n} \right) 
            \left( 1 + q^{3 + n} \right)} \times \\
& & q^{5 + 2n}\left( 1 - q^{1 + n} - q^{1 + 2n} + q^{2 + 2n} + q^{3 + 2n} - q^{6 + 2n}
            + q^{2 + 3n} + q^{7 + 3n} + q^{5 + 4n} + q^{6 + 4n} + 2q^{7 + 4n}
    \right.\\
& & \left.- q^{8 + 4n} - q^{9 + 4n} - q^{8 + 5n} + q^{9 + 6n} + q^{10 + 6n} 
     \right) J_{5_2}(1 + n)\\
& & + \frac {{\left( -1 + q^{2 + n} \right) }^2\left( 1 + q^{2 + n} \right)
              \left( -1 + q^{1 + 2n} \right) }
            {\left( -1 + q^n \right) \left( -1 + q^{3 + n} \right) 
             \left( 1 + q^{3 + n} \right) \left( -1 + q^{7 + 2n} \right)}
             \times  \\
& &  q\left( -1 + 2q^{2 + n} + q^{2 + 2n} - q^{4 + 2n} + q^{6 + 2n} 
            + q^{7 + 2n} - q^{4 + 3n} + q^{5 + 3n} + 3q^{6 + 3n} + 
            2q^{7 + 3n} - q^{8 + 3n} 
     \right.\\
& & - 2q^{9 + 3n} + q^{6 + 4n} - 2q^{8 + 4n} 
            - 2q^{9 + 4n} + q^{11 + 4n} - 
       q^{13 + 4n} - q^{8 + 5n} + 3q^{10 + 5n} + 3q^{11 + 5n}  
     \\
& &  - 3q^{13 + 5n}- 2  q^{14 + 5n} + q^{15 + 5n} + q^{10 + 6n} 
    + q^{11 + 6n} - q^{13 + 6n}      
     + 2q^{15 + 6n} + q^{16 + 6n} + q^{14 + 7n}  \\
& &  \left. + 3q^{15 + 7n}+ 2q^{16 + 7n} 
      - q^{17 + 7n} - q^{18 + 7n} - q^{17 + 8n} - q^{18 + 8n} + q^{20 + 9n} 
      \right) J_{5_2}(2 + n)\\
& & - \frac {\left( -1 + q^{1 + n} \right) \left( 1 + q^{1 + n} \right) 
            {\left( -1 + q^{3 + n} \right) }
             \left( -1 + q^{1 + 2n} \right)} 
            {\left( 1 + q^{n} \right)\left( -1 + q^{4 + n} \right) 
             \left( 1 + q^{4 + n} \right) \left( -1 + q^{7 + 2n} \right)}
             \times \\
& &  \left( -1 + q^{2 + n} + q^{3 + n} - q^{4 + 2n} - 3q^{5 + 2n} - 2q^{6 + 2n}             +q^{7 + 2n} + q^{8 + 2n} - q^{5 + 3n} - q^{6 + 3n} + q^{8 + 3n}
     \right.\\
& &  - 2q^{10 + 3n} - q^{11 + 3n} + q^{8 + 4n} - 3q^{10 + 4n} 
            - 3q^{11 + 4n} + 3q^{13 + 4n} + 2q^{14 + 4n} - q^{15 + 4n} 
            - q^{11 + 5n} \\
& &   + 2q^{13 + 5n} + 2q^{14 + 5n}
    - q^{16 + 5n} + q^{18 + 5n} + q^{14 + 6n} - q^{15 + 6n} 
            - 3q^{16 + 6n} - 2q^{17 + 6n} + q^{18 + 6n} \\
& &   \left.+ 2q^{19 + 6n}
            - q^{17 + 7n} + q^{19 + 7n} - q^{21 + 7n} - q^{22 + 7n} 
            - 2q^{22 + 8n} + q^{25 + 9n} \right) 
     J_{5_2}(3 + n)\\
& &  + \frac {\left( -1 + q^{1 + n} \right) 
     \left( 1 + q^{1 + n} \right) \left( -1 + q^{4 + n} \right) 
     \left( -1 + q^{1 + 2n} \right) 
     \left( -1 + q^{3 + 2n} \right)}
     {\left( -1 + q^n \right) \left( -1 + q^{3 + n} \right) 
     \left( 1 + q^{3 + n} \right)  \left( -1 + q^{7 + 2n} \right)
     \left( -1+q^{9 + 2n} \right)} \times \\
& &  q^{1 + n}
     \left( 1 + q - q^{4 + n} + q^{6 + 2n}
      + q^{7 + 2n} + 
       2q^{8 + 2n} - q^{9 + 2n} - q^{10 + 2n} + q^{8 + 3n} + q^{13 + 3n}
     \right.\\
& &   \left. - 
       q^{12 + 4n} + q^{13 + 4n} + q^{14 + 4n} - q^{17 + 4n} - q^{17 + 5n}
        + q^{21 + 6n}
       \right) J_{5_2}(4 + n)\\
& & + \frac {\left( -1 + q^{1 + n} \right)\left( 1 + q^{1 + n} \right) 
     \left( -1 + q^{2 + n} \right) \left( 1 + q^{2 + n} \right) 
     \left( -1 + q^{5 + n} \right) 
     \left( -1 + q^{1 + 2n} \right)\left( -1 + q^{3 + 2n} \right) }
     {\left( -1 + q^n \right) \left( -1 + q^{3 + n} \right) 
     \left( 1 + q^{3 + n} \right) \left( -1 + q^{4 + n} \right) 
     \left( 1 + q^{4 + n} \right) \left( -1 + q^{7 + 2n} \right)
     \left( -1+q^{9 + 2n} \right)}\times \\
& &     q^{4 + 2n} J_{5_2}(5 + n)
    = G(n), 
\end{eqnarray*}
where $G(n)$ is a sum of sums of proper $q$-hypergeometric function.
Setting $q=1$, and replacing $J_{5_2}(i + n)$ by $l^i$ and $q^n$ by $m^2$, 
the denominators cancell,  
the left hand side is changed to 
\begin{eqnarray*}
& & (1+m^2l)^2 \times \\
& & \{m^{14}+l(m^4-m^6+2 m^{10}+2 m^{12}-m^{14})
-l^2(-1+2m^2+2m^4-m^8+m^{10})+l^3\},
\end{eqnarray*}
and the second factor is equal to the $A$-polynomial of the knot $5_2$. 
So, this supports the  AJ conjecture for the knot $5_2$. 

%%%%%%%%%%%%%%%%%%%%
% 6_1 
%%%%%%%%%%%%%%%%%%%%%

 Furthermore, the colored Jones polynomial of the knot $6_1$ can be written 
\begin{equation}
J_{6_1}(n)=\sum_{k=0}^{\infty}\sum_{l=0}^k q^{-k^2-k}
       q^{nk}\frac {(q^{-1})_{n+k}(q^{-1})_{n-1}}
                    {(q^{-1})_{n}(q^{-1})_{n-k-1}}
       q^{l(k+1)}\frac {(q)_k}{(q)_{l} (q)_{k-l}}. 
\end{equation}
Using {\tt qMultiSum.m} again, we can obtain 
the inhomogeneous recursive relation of $J_{6_1}(n)$ 
\begin{eqnarray*}
& &q^{8 + 4n} J_{6_1}(n) \\
& & - \frac {{\left( -1 + q^{1 + n} \right) }^2\left( 1 + q^{1 + n} \right) }
            {\left( -1 + q^n \right) \left( -1 + q^{3 + n} \right) 
             \left( 1 + q^{3 + n} \right)}
      \times \\
& & q^{5 + 2n}\left( 1 + q - 2q^{2 + n} - q^{3 + n} - q^{1 + 2n} - 
       q^{2 + 2n} + q^{3 + 2n} - q^{6 + 2n} - q^{7 + 2n} + q^{3 + 3n}
        - q^{5 + 3n}\right. \\
& &   \left. - q^{6 + 3n} - q^{7 + 3n} + q^{8 + 3n} - q^{4 + 4n} 
         + q^{7 + 4n} + q^{8 + 4n} - 
       q^{9 + 4n} - q^{9 + 5n} + q^{10 + 6n} \right) 
       J_{6_1}(1 + n) \\
& & - \frac {{\left( -1 + q^{2 + n} \right) }^2\left( 1 + q^{2 + n} \right)
             \left( -1 + q^{1 + 2n} \right)}
            {\left( -1 + q^n \right) \left( -1 + q^{3 + n} \right) 
             \left( 1 + q^{3 + n} \right)\left( -1 + q^{7 + 2n} \right) }
       \times \\
& &   q \left( -1 + 2q^{2 + n} + q^{3 + n} + q^{2 + 2n} + q^{3 + 2n} 
      - q^{4 + 2n} - 2q^{5 + 2n} + q^{6 + 2n} + q^{7 + 2n} - 2q^{4 + 3n}
       - q^{5 + 3n} 
      \right.\\
& &  + 2q^{6 + 3n} + 3q^{7 + 3n} - 2q^{9 + 3n} - q^{10 + 3n} - q^{5 + 4n} 
     + q^{6 + 4n} + q^{7 + 4n} - 2q^{8 + 4n} - 5q^{9 + 4n} - 3q^{10 + 4n} 
     \\
& &  + q^{12 + 4n} - q^{13 + 4n} - 2q^{8 + 5n} - 3q^{9 + 5n} + 3q^{11 + 5n} 
     + q^{12 + 5n} 
     - 2q^{13 + 5n} - 3q^{14 + 5n} + 2q^{11 + 6n}\\ 
& &  + q^{12 + 6n}- 2q^{13 + 6n} - 2q^{14 + 6n} 
     + 2q^{16 + 6n} - q^{13 + 7n} + 2q^{15 + 7n} + 3q^{16 + 7n} - q^{18 + 7n}\\
& &    \left. - 2q^{18 + 8n} + q^{20 + 9n} \right) 
      J_{6_1}(2 + n) \\
& & + \frac {\left( -1 + q^{1 + n} \right) \left( 1 + q^{1 + n} \right) 
             \left( -1 + q^{3 + n} \right) \left( -1 + q^{1 + 2n} \right) }
            {\left( -1 + q^n \right) \left( -1 + q^{4 + n} \right) 
             \left( 1 + q^{4 + n} \right) \left( -1 + q^{7 + 2n} \right) } 
      \times \\
& &   \left( -1 + 2q^{3 + n} + q^{3 + 2n} - 2q^{5 + 2n} - 3q^{6 + 2n} 
       + q^{8 + 2n} - 2q^{6 + 3n} - q^{7 + 3n} + 2q^{8 + 3n} + 2q^{9 + 3n} 
       \right.\\
& &   - 2q^{11 + 3n}+ 2q^{8 + 4n} + 3q^{9 + 4n} - 3q^{11 + 4n} - q^{12 + 4n} 
      + 2q^{13 + 4n}
      + 3q^{14 + 4n} + q^{10 + 5n} - q^{11 + 5n}   \\
&&    - q^{12 + 5n}+ 2q^{13 + 5n} + 5q^{14 + 5n}+ 3q^{15 + 5n} - q^{17 + 5n}
       + q^{18 + 5n}  + 2q^{14 + 6n} + q^{15 + 6n} - 2q^{16 + 6n} \\
& &  - 3q^{17 + 6n}+ 2q^{19 + 6n} + q^{20 + 6n} - q^{17 + 7n} - q^{18 + 7n}
     + q^{19 + 7n} + 2q^{20 + 7n} - q^{21 + 7n} - q^{22 + 7n}  
     \\
& &  \left. - 2q^{22 + 8n}- q^{23 + 8n} + q^{25 + 9n} \right) 
     J_{6_1}(3 + n) \\
& & - \frac {\left( -1 + q^{1 + n} \right) \left( 1 + q^{1 + n} \right) 
             \left( -1 + q^{4 + n} \right)\left( -1 + q^{1 + 2n} \right) 
             \left( -1 + q^{3 + 2n} \right) }
            {\left( -1 + q^n \right) \left( -1 + q^{3 + n} \right) 
             \left( 1 + q^{3 + n} \right) \left( -1 + q^{7 + 2n} \right) 
            \left( -1 + q^{9 + 2n} \right) } \times \\
& &  q^{2 + n}
     \left( 1 - q^{4 + n} - q^{4 + 2n} + q^{7 + 2n} + q^{8 + 2n} 
     - q^{9 + 2n} + q^{8 + 3n} - q^{10 + 3n} - q^{11 + 3n} - q^{12 + 3n} 
     \right.\\
& &   + q^{13 + 3n} - q^{11 + 4n} - q^{12 + 4n} + q^{13 + 4n} - q^{16 + 4n} 
     - q^{17 + 4n} - 2q^{17 + 5n} - q^{18 + 5n}
     \\
& & \left.  + q^{20 + 6n} + q^{21 + 6n} \right) 
    J_{6_1}(4 + n) \\
& & + \frac {\left( -1 + q^{1 + n} \right) \left( 1 + q^{1 + n} \right) 
             \left( -1 + q^{2 + n} \right) \left( 1 + q^{2 + n} \right) 
             \left( -1 + q^{5 + n} \right) \left( -1 + q^{1 + 2n} \right)
             \left( -1 + q^{3 + 2n} \right)}
            {\left( -1 + q^n \right) \left( -1 + q^{3 + n} \right) 
             \left( 1 + q^{3 + n} \right) \left( -1 + q^{4 + n} \right)
             \left( 1 + q^{4 + n} \right) \left( -1 + q^{7 + 2n} \right) 
              \left( -1 + q^{9 + 2n} \right)} \times \\
& &     q^{18 + 5n}J_{6_1}(5 + n) = G'(n), 
\end{eqnarray*}
where $G'(n)$ is a sum of sums of proper $q$-hypergeometric function.
 Setting $q=1$, and replacing $J_{6_1}(i + n)$ by $l^i$ and $q^n$ by $m^2$, 
 the left hand side is changed to 
\begin{eqnarray*}
& & (1+m^2l) 
    \{ m^{8}+l(-2m^4+3m^6+3m^8+2 m^{14}-m^{16})\\
& &     +l^2(1-3m^2-m^4+3m^6+6m^8+3m^{10}-m^{12}-3m^{14}+m^{16})\\
& &     +l^3(-1+m^2+m^4+3m^{10}+3m^{12}-2m^{14})+m^8l^4\},
\end{eqnarray*}
and the second factor coincides with the A-polynomial of the knot $6_1$.

%%%%%%%%%%%%%%%%%%%%%%%%%%%%%%%%%%%%
%volume conjecture
%%%%%%%%%%%%%%%%%%%%%%%%%%%%%%%%%%%%

\quad 

 Finally, we discuss our result in terms of 
{\it  Volume conjecture} 
due to Kashaev, 
H.Murakami and J.Murakami. 
(see \cite{Kashaev:LETMP97},\cite{Murakami/Murakami:volume}). 
Let $K$ be a hyperbolic knot in $S^3$ and $M_K$ the complement of $K$.
We write  $\hat{J}_K(n)$ for $J_K(n)$ evaluated at $q=\exp\frac{ 2\pi \I} n$.
Then, {\it  Volume conjecture} is 
\begin{conj}
\[
2\pi\lim_{n\to \infty} \frac {\log|\hat{J}_K(n)|}{n}=Vol(M_K),
\]
where $Vol(M_K)$ is the hyperbolic volume of $M_K$.
\end{conj}
From the formula of  $J_{K_p}(n)$  given by Masbaum, we have  
\begin{equation}
\hat{J}_{K_p}(n)=\sum_{k=0}^{\infty}\sum_{l=0}^k (-1)^{k+l} q^{k(k+3)/2}
        (q)_k(q^{-1})_k
       q^{l(l+1)p+l(l-1)/2}\frac {(q)_{2l+1}}  {(q)_{2l}}
       \frac {(q)_k}{(q)_{k+l+1} (q)_{k-l}}.
\end{equation}
Computing the optimistic limit \cite{HMurakami:RIMS00} of $\hat{J}_K(n)$, 
it is conjectured that 
there exists a solution $(x_0,y_0)$ 
to the two equations 
\begin{equation} 
f_1(1,1,x,y)=1, \; f_2(1,1,x,y)=1
\label{eqn:critical}
\end{equation}
such that 
\[
Vol(M_K)=3D(x_0)-D(x_0y_0)-D(\frac {x_0} {y_0}), 
\] 
where $f_i$ is the same notation as in Section 2, and 
$D(z)=\mathrm{Im }(\li_2(z))+\log |z| \arg (1-z)$ with 
$\li_2(z)=-\int_{0}^z \frac {\log (1-w)} w dw$. 
By using Mathematica 4 and Maple V, we can numerically see    
that for $p=-1, 2,-2$, this conjecture is true. 

In \cite{Yokota:volume}, Y.Yokota introduced a way 
to construct an ideal triangulation of the complement $M_K$ of $K$ associated to  Kashaev's $R$-matrices, 
and  gave a sketch proof of  {\it Volume conjecture}. 
Moreover, in \cite{Yokota:IIS03}, he explained a method to calculate  
the $A$-polynomial of $K$ from the triangulation of $\partial M_K$
 in \cite{Yokota:volume}. 
We hope to find a geometrical interpretation of 
 our method to obtain the $A$-polynomial, comparing with Yokota's method.

%%%%%%%%%%%%%%%%%%%%%%%%%%%%%%%%%%%%%%%%%%%%%%%%%%%%%%%%%%%%%%%%%%%%%
\bibliography{mrabbrev,toshie-bib}
\bibliographystyle{amsplain}
%%%%%%%%%%%%%%%%%%%%%%%%%%%%%%%%%%%%%%%%%%%%%%%%%%%%%%%%%%%%%%%%%%%%%

\end{document}